\documentclass{gtart}

\def\ifplaintex{\expandafter\ifx\csname documentclass\endcsname\relax}

\def\gtp{{\mathsurround=0pt\it $\cal G\mskip-2mu$eometry \&\ 
$\cal T\!\!$opology $\cal P\!$ublications}}  

\def\Addressesr{\bigskip
{\small \parskip 0pt \leftskip 0pt \rightskip 0pt plus 1fil \def\\{\par}
\sl\theaddress\par
\medskip
\rm Email:\stdspace\tt\theemail\hfill\rm Received:\qua\receiveddate \par}}

\def\recd{{\small Received:\qua\receiveddate\ifx\reviseddate\relax
\else\qquad Revised:\qua\reviseddate\fi\par}} 

\def\lognumber#1{\def\thelognumber{#1}}
\def\volumenumber#1{\def\thevolumenumber{#1}}
\def\volumeyear#1{\def\thevolumeyear{#1}}
\def\papernumber#1{\def\thepapernumber{#1}}
\def\pagenumbers#1#2{\def\startpage{#1}\def\finishpage{#2}}
\def\published#1{\def\publishdate{#1}}

\def\received#1{\def\receiveddate{#1}}

\def\accepted#1{\def\accepteddate{#1}}

\long\def\asciiabstract#1{\long\def\theasciiabstract{#1}}

\let\\\par\let\thelognumber\relax\let\thevolumenumber\relax
\let\thepapernumber\relax\let\thevolumeyear\relax\let\startpage\relax
\let\finishpage\relax\let\publishdate\relax\let\receiveddate\relax
\let\reviseddate\relax\let\accepteddate\relax\let\theasciititle\relax
\let\theasciiauthors\relax
\let\theasciiabstract\relax

\let\theasciiemail\relax

\ifplaintex
\font\logobig=cmssbx10 scaled 3836
\font\logomed=cmssbx10 scaled 2557
\else
\font\logobig=cmssbx10 scaled 4200
\font\logomed=cmssbx10 scaled 2800
\fi

\long\def\makeagttitle{   
\count0=\startpage
\agt\hfill      
\hbox to 45truept{\vbox to 0pt{\vglue -13truept{\logomed A\kern -.37em{\logobig 
T}\kern -.38em G}\vss}\hss}
\break
{\small Volume \thevolumenumber\ (\thevolumeyear)
\startpage--\finishpage\nl
Published: \publishdate}

\vglue .25truein

{\parskip=0pt\leftskip 0pt plus
1fil\def\\{\par\smallskip}{\Large\bf\thetitle}\par\medskip} \vglue
0.05truein

{\parskip=0pt\leftskip 0pt plus 1fil\def\\{\par}{\sc\theauthors}
\par\medskip}%
 
\vglue 0.03truein 

{\small\leftskip 25truept\rightskip 25truept{\bf Abstract}\stdspace\theabstract

{\bf AMS Classification}\stdspace\theprimaryclass
\ifx\thesecondaryclass\relax\else; \thesecondaryclass\fi\par
{\bf Keywords}\stdspace \thekeywords\par}\vglue 7truept

}   

\ifplaintex
\font\phead=cmsl9 scaled 950
\font\pnum=cmbx10 scaled 913
\font\pfoot=cmsl9 scaled 950
\headline{\vbox to 0pt{\vskip -4.5mm\line{\small\phead\ifnum
\count0=\startpage ISSN 1472-2739 (on-line) 1472-2747 (printed)
\hfill {\pnum\folio}\else\ifodd\count0\def\\{ }%
\ifx\theshorttitle\relax\thetitle\else\theshorttitle\fi\hfill{\pnum\folio}
\else\def\\{ and }{\pnum\folio}\hfill\ifx\theshortauthors\relax\theauthors
\else\theshortauthors\fi\fi\fi}\vss}}
\footline{\vbox to 0pt{\vglue 0mm\line{\small\pfoot\ifnum\count0=\startpage
\copyright\ \gtp\hfill\else
\agt, Volume \thevolumenumber\ (\thevolumeyear)\hfill\fi}\vss}}
\else
\font=cmsl9 scaled 1050
\font\lnum=cmbx10 
\font=cmsl9 scaled 1050
\makeatletter
\def\@oddhead{{\small\lhead\ifnum\count0=\startpage ISSN 1472-2739 
(on-line) 1472-2747 (printed)\hfill {\lnum\number\count0}\else\ifodd\count0
\def\\{ }\ifx\theshorttitle\relax \thetitle \else\theshorttitle\fi\hfill
{\lnum\number\count0}\else\def\\{ and }{\lnum\number\count0}
\hfill\ifx\theshortauthors\relax 
\theauthors\else\theshortauthors\fi\fi\fi}}\def\@evenhead{\@oddhead}
\def\@oddfoot{\small\lfoot\ifnum\count0=\startpage\copyright\ \gtp\hfill\else
\agt, Volume \thevolumenumber\ (\thevolumeyear)\hfill\fi}
\def\@evenfoot{\@oddfoot}
\makeatother
\fi
\let\maketitlepage\makeagttitle
\let\makeshorttitle\maketitlepage
\let\maketitle\maketitlepage

\newwrite\gtoutfile
\long\gdef\makeheadfile{  
{\def\\{, }\def\s{ }
\immediate\openout\gtoutfile head.xxx
\immediate\write\gtoutfile{To: math@arxiv.org}
\immediate\write\gtoutfile{Subject: put OR rep NNNNN:ppppp}
\immediate\write\gtoutfile{--text follows this line--}
\immediate\write\gtoutfile{Proxy-for: \ifx\theasciiauthors\relax
\theauthors\else\theasciiauthors\fi\s<\ifx\theasciiemail\relax\theemail\else\theasciiemail\fi>}
\immediate\write\gtoutfile{\noexpand\\}
\immediate\write\gtoutfile{Authors: \ifx\theasciiauthors\relax
\theauthors\else\theasciiauthors\fi}
{\def\\{ }\immediate\write\gtoutfile{Title: \ifx\theasciititle\relax
\thetitle\else\theasciititle\fi}}
\immediate\write\gtoutfile{Subj-class: GT or SG, GR etc}
\immediate\write\gtoutfile{MSC-class: \theprimaryclass\ifx\thesecondaryclass\relax\else, \thesecondaryclass\fi}
\immediate\write\gtoutfile{Journal-ref: Algebr. Geom. Topol. \thevolumenumber\s
(\thevolumeyear) \startpage-\finishpage}
\immediate\write\gtoutfile{Comments: Published by Algebraic and
Geometric Topology at}
\immediate\write\gtoutfile{\s\s\s  http://www.maths.warwick.ac.uk/agt/AGTVol\thevolumenumber/agt-\thevolumenumber-\thepapernumber.abs.html}
\immediate\write\gtoutfile{\noexpand\\}
\immediate\write\gtoutfile{}
\ifx\theasciiabstract\relax
\immediate\write\gtoutfile{\theabstract}\else
\immediate\write\gtoutfile{\theasciiabstract}\fi
\immediate\write\gtoutfile{}
\immediate\write\gtoutfile{\noexpand\\}
\immediate\write\gtoutfile{}
\immediate\closeout\gtoutfile}}  

\def\maketitlepage{\makeagttitle\makeheadfile}
\let\makeshorttitle\maketitlepage
\let\maketitle\maketitlepage

\def\ifplaintex{\expandafter\ifx\csname documentclass\endcsname\relax}

\def\gtp{{\mathsurround=0pt\it $\cal G\mskip-2mu$eometry \&\ 
$\cal T\!\!$opology $\cal P\!$ublications}}  

\def\Addressesr{\bigskip
{\small \parskip 0pt \leftskip 0pt \rightskip 0pt plus 1fil \def\\{\par}
\sl\theaddress\par
\medskip
\rm Email:\stdspace\tt\theemail\hfill\rm Received:\qua\receiveddate \par}}

\def\recd{{\small Received:\qua\receiveddate\ifx\reviseddate\relax
\else\qquad Revised:\qua\reviseddate\fi\par}} 

\def\lognumber#1{\def\thelognumber{#1}}
\def\volumenumber#1{\def\thevolumenumber{#1}}
\def\volumeyear#1{\def\thevolumeyear{#1}}
\def\papernumber#1{\def\thepapernumber{#1}}
\def\pagenumbers#1#2{\def\startpage{#1}\def\finishpage{#2}}
\def\published#1{\def\publishdate{#1}}

\def\received#1{\def\receiveddate{#1}}

\def\accepted#1{\def\accepteddate{#1}}

\long\def\asciiabstract#1{\long\def\theasciiabstract{#1}}

\let\\\par\let\thelognumber\relax\let\thevolumenumber\relax
\let\thepapernumber\relax\let\thevolumeyear\relax\let\startpage\relax
\let\finishpage\relax\let\publishdate\relax\let\receiveddate\relax
\let\reviseddate\relax\let\accepteddate\relax\let\theasciititle\relax
\let\theasciiauthors\relax
\let\theasciiabstract\relax

\let\theasciiemail\relax

\ifplaintex
\font\logobig=cmssbx10 scaled 3836
\font\logomed=cmssbx10 scaled 2557
\else
\font\logobig=cmssbx10 scaled 4200
\font\logomed=cmssbx10 scaled 2800
\fi

\long\def\makeagttitle{   
\count0=\startpage
\agt\hfill      
\hbox to 45truept{\vbox to 0pt{\vglue -13truept{\logomed A\kern -.37em{\logobig 
T}\kern -.38em G}\vss}\hss}
\break
{\small Volume \thevolumenumber\ (\thevolumeyear)
\startpage--\finishpage\nl
Published: \publishdate}

\vglue .25truein

{\parskip=0pt\leftskip 0pt plus
1fil\def\\{\par\smallskip}{\Large\bf\thetitle}\par\medskip} \vglue
0.05truein

{\parskip=0pt\leftskip 0pt plus 1fil\def\\{\par}{\sc\theauthors}
\par\medskip}%
 
\vglue 0.03truein 

{\small\leftskip 25truept\rightskip 25truept{\bf Abstract}\stdspace\theabstract

{\bf AMS Classification}\stdspace\theprimaryclass
\ifx\thesecondaryclass\relax\else; \thesecondaryclass\fi\par
{\bf Keywords}\stdspace \thekeywords\par}\vglue 7truept

}   

\ifplaintex
\font\phead=cmsl9 scaled 950
\font\pnum=cmbx10 scaled 913
\font\pfoot=cmsl9 scaled 950
\headline{\vbox to 0pt{\vskip -4.5mm\line{\small\phead\ifnum
\count0=\startpage ISSN 1472-2739 (on-line) 1472-2747 (printed)
\hfill {\pnum\folio}\else\ifodd\count0\def\\{ }%
\ifx\theshorttitle\relax\thetitle\else\theshorttitle\fi\hfill{\pnum\folio}
\else\def\\{ and }{\pnum\folio}\hfill\ifx\theshortauthors\relax\theauthors
\else\theshortauthors\fi\fi\fi}\vss}}
\footline{\vbox to 0pt{\vglue 0mm\line{\small\pfoot\ifnum\count0=\startpage
\copyright\ \gtp\hfill\else
\agt, Volume \thevolumenumber\ (\thevolumeyear)\hfill\fi}\vss}}
\else
\font=cmsl9 scaled 1050
\font\lnum=cmbx10 
\font=cmsl9 scaled 1050
\makeatletter
\def\@oddhead{{\small\lhead\ifnum\count0=\startpage ISSN 1472-2739 
(on-line) 1472-2747 (printed)\hfill {\lnum\number\count0}\else\ifodd\count0
\def\\{ }\ifx\theshorttitle\relax \thetitle \else\theshorttitle\fi\hfill
{\lnum\number\count0}\else\def\\{ and }{\lnum\number\count0}
\hfill\ifx\theshortauthors\relax 
\theauthors\else\theshortauthors\fi\fi\fi}}\def\@evenhead{\@oddhead}
\def\@oddfoot{\small\lfoot\ifnum\count0=\startpage\copyright\ \gtp\hfill\else
\agt, Volume \thevolumenumber\ (\thevolumeyear)\hfill\fi}
\def\@evenfoot{\@oddfoot}
\makeatother
\fi
\let\maketitlepage\makeagttitle
\let\makeshorttitle\maketitlepage
\let\maketitle\maketitlepage

\newwrite\gtoutfile
\long\gdef\makeheadfile{  
{\def\\{, }\def\s{ }
\immediate\openout\gtoutfile head.xxx
\immediate\write\gtoutfile{To: math@arxiv.org}
\immediate\write\gtoutfile{Subject: put OR rep NNNNN:ppppp}
\immediate\write\gtoutfile{--text follows this line--}
\immediate\write\gtoutfile{Proxy-for: \ifx\theasciiauthors\relax
\theauthors\else\theasciiauthors\fi\s<\ifx\theasciiemail\relax\theemail\else\theasciiemail\fi>}
\immediate\write\gtoutfile{\noexpand\\}
\immediate\write\gtoutfile{Authors: \ifx\theasciiauthors\relax
\theauthors\else\theasciiauthors\fi}
{\def\\{ }\immediate\write\gtoutfile{Title: \ifx\theasciititle\relax
\thetitle\else\theasciititle\fi}}
\immediate\write\gtoutfile{Subj-class: GT or SG, GR etc}
\immediate\write\gtoutfile{MSC-class: \theprimaryclass\ifx\thesecondaryclass\relax\else, \thesecondaryclass\fi}
\immediate\write\gtoutfile{Journal-ref: Algebr. Geom. Topol. \thevolumenumber\s
(\thevolumeyear) \startpage-\finishpage}
\immediate\write\gtoutfile{Comments: Published by Algebraic and
Geometric Topology at}
\immediate\write\gtoutfile{\s\s\s  http://www.maths.warwick.ac.uk/agt/AGTVol\thevolumenumber/agt-\thevolumenumber-\thepapernumber.abs.html}
\immediate\write\gtoutfile{\noexpand\\}
\immediate\write\gtoutfile{}
\ifx\theasciiabstract\relax
\immediate\write\gtoutfile{\theabstract}\else
\immediate\write\gtoutfile{\theasciiabstract}\fi
\immediate\write\gtoutfile{}
\immediate\write\gtoutfile{\noexpand\\}
\immediate\write\gtoutfile{}
\immediate\closeout\gtoutfile}}  

\def\maketitlepage{\makeagttitle\makeheadfile}
\let\makeshorttitle\maketitlepage
\let\maketitle\maketitlepage

\def\ifplaintex{\expandafter\ifx\csname documentclass\endcsname\relax}

\def\gtp{{\mathsurround=0pt\it $\cal G\mskip-2mu$eometry \&\ 
$\cal T\!\!$opology $\cal P\!$ublications}}  

\def\Addressesr{\bigskip
{\small \parskip 0pt \leftskip 0pt \rightskip 0pt plus 1fil \def\\{\par}
\sl\theaddress\par
\medskip
\rm Email:\stdspace\tt\theemail\hfill\rm Received:\qua\receiveddate \par}}

\def\recd{{\small Received:\qua\receiveddate\ifx\reviseddate\relax
\else\qquad Revised:\qua\reviseddate\fi\par}} 

\def\lognumber#1{\def\thelognumber{#1}}
\def\volumenumber#1{\def\thevolumenumber{#1}}
\def\volumeyear#1{\def\thevolumeyear{#1}}
\def\papernumber#1{\def\thepapernumber{#1}}
\def\pagenumbers#1#2{\def\startpage{#1}\def\finishpage{#2}}
\def\published#1{\def\publishdate{#1}}

\def\received#1{\def\receiveddate{#1}}

\def\accepted#1{\def\accepteddate{#1}}

\long\def\asciiabstract#1{\long\def\theasciiabstract{#1}}

\let\\\par\let\thelognumber\relax\let\thevolumenumber\relax
\let\thepapernumber\relax\let\thevolumeyear\relax\let\startpage\relax
\let\finishpage\relax\let\publishdate\relax\let\receiveddate\relax
\let\reviseddate\relax\let\accepteddate\relax\let\theasciititle\relax
\let\theasciiauthors\relax
\let\theasciiabstract\relax

\let\theasciiemail\relax

\ifplaintex
\font\logobig=cmssbx10 scaled 3836
\font\logomed=cmssbx10 scaled 2557
\else
\font\logobig=cmssbx10 scaled 4200
\font\logomed=cmssbx10 scaled 2800
\fi

\long\def\makeagttitle{   
\count0=\startpage
\agt\hfill      
\hbox to 45truept{\vbox to 0pt{\vglue -13truept{\logomed A\kern -.37em{\logobig 
T}\kern -.38em G}\vss}\hss}
\break
{\small Volume \thevolumenumber\ (\thevolumeyear)
\startpage--\finishpage\nl
Published: \publishdate}

\vglue .25truein

{\parskip=0pt\leftskip 0pt plus
1fil\def\\{\par\smallskip}{\Large\bf\thetitle}\par\medskip} \vglue
0.05truein

{\parskip=0pt\leftskip 0pt plus 1fil\def\\{\par}{\sc\theauthors}
\par\medskip}%
 
\vglue 0.03truein 

{\small\leftskip 25truept\rightskip 25truept{\bf Abstract}\stdspace\theabstract

{\bf AMS Classification}\stdspace\theprimaryclass
\ifx\thesecondaryclass\relax\else; \thesecondaryclass\fi\par
{\bf Keywords}\stdspace \thekeywords\par}\vglue 7truept

}   

\ifplaintex
\font\phead=cmsl9 scaled 950
\font\pnum=cmbx10 scaled 913
\font\pfoot=cmsl9 scaled 950
\headline{\vbox to 0pt{\vskip -4.5mm\line{\small\phead\ifnum
\count0=\startpage ISSN 1472-2739 (on-line) 1472-2747 (printed)
\hfill {\pnum\folio}\else\ifodd\count0\def\\{ }%
\ifx\theshorttitle\relax\thetitle\else\theshorttitle\fi\hfill{\pnum\folio}
\else\def\\{ and }{\pnum\folio}\hfill\ifx\theshortauthors\relax\theauthors
\else\theshortauthors\fi\fi\fi}\vss}}
\footline{\vbox to 0pt{\vglue 0mm\line{\small\pfoot\ifnum\count0=\startpage
\copyright\ \gtp\hfill\else
\agt, Volume \thevolumenumber\ (\thevolumeyear)\hfill\fi}\vss}}
\else
\font=cmsl9 scaled 1050
\font\lnum=cmbx10 
\font=cmsl9 scaled 1050
\makeatletter
\def\@oddhead{{\small\lhead\ifnum\count0=\startpage ISSN 1472-2739 
(on-line) 1472-2747 (printed)\hfill {\lnum\number\count0}\else\ifodd\count0
\def\\{ }\ifx\theshorttitle\relax \thetitle \else\theshorttitle\fi\hfill
{\lnum\number\count0}\else\def\\{ and }{\lnum\number\count0}
\hfill\ifx\theshortauthors\relax 
\theauthors\else\theshortauthors\fi\fi\fi}}\def\@evenhead{\@oddhead}
\def\@oddfoot{\small\lfoot\ifnum\count0=\startpage\copyright\ \gtp\hfill\else
\agt, Volume \thevolumenumber\ (\thevolumeyear)\hfill\fi}
\def\@evenfoot{\@oddfoot}
\makeatother
\fi
\let\maketitlepage\makeagttitle
\let\makeshorttitle\maketitlepage
\let\maketitle\maketitlepage

\newwrite\gtoutfile
\long\gdef\makeheadfile{  
{\def\\{, }\def\s{ }
\immediate\openout\gtoutfile head.xxx
\immediate\write\gtoutfile{To: math@arxiv.org}
\immediate\write\gtoutfile{Subject: put OR rep NNNNN:ppppp}
\immediate\write\gtoutfile{--text follows this line--}
\immediate\write\gtoutfile{Proxy-for: \ifx\theasciiauthors\relax
\theauthors\else\theasciiauthors\fi\s<\ifx\theasciiemail\relax\theemail\else\theasciiemail\fi>}
\immediate\write\gtoutfile{\noexpand\\}
\immediate\write\gtoutfile{Authors: \ifx\theasciiauthors\relax
\theauthors\else\theasciiauthors\fi}
{\def\\{ }\immediate\write\gtoutfile{Title: \ifx\theasciititle\relax
\thetitle\else\theasciititle\fi}}
\immediate\write\gtoutfile{Subj-class: GT or SG, GR etc}
\immediate\write\gtoutfile{MSC-class: \theprimaryclass\ifx\thesecondaryclass\relax\else, \thesecondaryclass\fi}
\immediate\write\gtoutfile{Journal-ref: Algebr. Geom. Topol. \thevolumenumber\s
(\thevolumeyear) \startpage-\finishpage}
\immediate\write\gtoutfile{Comments: Published by Algebraic and
Geometric Topology at}
\immediate\write\gtoutfile{\s\s\s  http://www.maths.warwick.ac.uk/agt/AGTVol\thevolumenumber/agt-\thevolumenumber-\thepapernumber.abs.html}
\immediate\write\gtoutfile{\noexpand\\}
\immediate\write\gtoutfile{}
\ifx\theasciiabstract\relax
\immediate\write\gtoutfile{\theabstract}\else
\immediate\write\gtoutfile{\theasciiabstract}\fi
\immediate\write\gtoutfile{}
\immediate\write\gtoutfile{\noexpand\\}
\immediate\write\gtoutfile{}
\immediate\closeout\gtoutfile}}  

\def\maketitlepage{\makeagttitle\makeheadfile}
\let\makeshorttitle\maketitlepage
\let\maketitle\maketitlepage

\def\ifplaintex{\expandafter\ifx\csname documentclass\endcsname\relax}

\def\gtp{{\mathsurround=0pt\it $\cal G\mskip-2mu$eometry \&\ 
$\cal T\!\!$opology $\cal P\!$ublications}}  

\def\Addressesr{\bigskip
{\small \parskip 0pt \leftskip 0pt \rightskip 0pt plus 1fil \def\\{\par}
\sl\theaddress\par
\medskip
\rm Email:\stdspace\tt\theemail\hfill\rm Received:\qua\receiveddate \par}}

\def\recd{{\small Received:\qua\receiveddate\ifx\reviseddate\relax
\else\qquad Revised:\qua\reviseddate\fi\par}} 

\def\lognumber#1{\def\thelognumber{#1}}
\def\volumenumber#1{\def\thevolumenumber{#1}}
\def\volumeyear#1{\def\thevolumeyear{#1}}
\def\papernumber#1{\def\thepapernumber{#1}}
\def\pagenumbers#1#2{\def\startpage{#1}\def\finishpage{#2}}
\def\published#1{\def\publishdate{#1}}

\def\received#1{\def\receiveddate{#1}}

\def\accepted#1{\def\accepteddate{#1}}

\long\def\asciiabstract#1{\long\def\theasciiabstract{#1}}

\let\\\par\let\thelognumber\relax\let\thevolumenumber\relax
\let\thepapernumber\relax\let\thevolumeyear\relax\let\startpage\relax
\let\finishpage\relax\let\publishdate\relax\let\receiveddate\relax
\let\reviseddate\relax\let\accepteddate\relax\let\theasciititle\relax
\let\theasciiauthors\relax
\let\theasciiabstract\relax

\let\theasciiemail\relax

\ifplaintex
\font\logobig=cmssbx10 scaled 3836
\font\logomed=cmssbx10 scaled 2557
\else
\font\logobig=cmssbx10 scaled 4200
\font\logomed=cmssbx10 scaled 2800
\fi

\long\def\makeagttitle{   
\count0=\startpage
\agt\hfill      
\hbox to 45truept{\vbox to 0pt{\vglue -13truept{\logomed A\kern -.37em{\logobig 
T}\kern -.38em G}\vss}\hss}
\break
{\small Volume \thevolumenumber\ (\thevolumeyear)
\startpage--\finishpage\nl
Published: \publishdate}

\vglue .25truein

{\parskip=0pt\leftskip 0pt plus
1fil\def\\{\par\smallskip}{\Large\bf\thetitle}\par\medskip} \vglue
0.05truein

{\parskip=0pt\leftskip 0pt plus 1fil\def\\{\par}{\sc\theauthors}
\par\medskip}%
 
\vglue 0.03truein 

{\small\leftskip 25truept\rightskip 25truept{\bf Abstract}\stdspace\theabstract

{\bf AMS Classification}\stdspace\theprimaryclass
\ifx\thesecondaryclass\relax\else; \thesecondaryclass\fi\par
{\bf Keywords}\stdspace \thekeywords\par}\vglue 7truept

}   

\ifplaintex
\font\phead=cmsl9 scaled 950
\font\pnum=cmbx10 scaled 913
\font\pfoot=cmsl9 scaled 950
\headline{\vbox to 0pt{\vskip -4.5mm\line{\small\phead\ifnum
\count0=\startpage ISSN 1472-2739 (on-line) 1472-2747 (printed)
\hfill {\pnum\folio}\else\ifodd\count0\def\\{ }%
\ifx\theshorttitle\relax\thetitle\else\theshorttitle\fi\hfill{\pnum\folio}
\else\def\\{ and }{\pnum\folio}\hfill\ifx\theshortauthors\relax\theauthors
\else\theshortauthors\fi\fi\fi}\vss}}
\footline{\vbox to 0pt{\vglue 0mm\line{\small\pfoot\ifnum\count0=\startpage
\copyright\ \gtp\hfill\else
\agt, Volume \thevolumenumber\ (\thevolumeyear)\hfill\fi}\vss}}
\else
\font=cmsl9 scaled 1050
\font\lnum=cmbx10 
\font=cmsl9 scaled 1050
\makeatletter
\def\@oddhead{{\small\lhead\ifnum\count0=\startpage ISSN 1472-2739 
(on-line) 1472-2747 (printed)\hfill {\lnum\number\count0}\else\ifodd\count0
\def\\{ }\ifx\theshorttitle\relax \thetitle \else\theshorttitle\fi\hfill
{\lnum\number\count0}\else\def\\{ and }{\lnum\number\count0}
\hfill\ifx\theshortauthors\relax 
\theauthors\else\theshortauthors\fi\fi\fi}}\def\@evenhead{\@oddhead}
\def\@oddfoot{\small\lfoot\ifnum\count0=\startpage\copyright\ \gtp\hfill\else
\agt, Volume \thevolumenumber\ (\thevolumeyear)\hfill\fi}
\def\@evenfoot{\@oddfoot}
\makeatother
\fi
\let\maketitlepage\makeagttitle
\let\makeshorttitle\maketitlepage
\let\maketitle\maketitlepage

\newwrite\gtoutfile
\long\gdef\makeheadfile{  
{\def\\{, }\def\s{ }
\immediate\openout\gtoutfile head.xxx
\immediate\write\gtoutfile{To: math@arxiv.org}
\immediate\write\gtoutfile{Subject: put OR rep NNNNN:ppppp}
\immediate\write\gtoutfile{--text follows this line--}
\immediate\write\gtoutfile{Proxy-for: \ifx\theasciiauthors\relax
\theauthors\else\theasciiauthors\fi\s<\ifx\theasciiemail\relax\theemail\else\theasciiemail\fi>}
\immediate\write\gtoutfile{\noexpand\\}
\immediate\write\gtoutfile{Authors: \ifx\theasciiauthors\relax
\theauthors\else\theasciiauthors\fi}
{\def\\{ }\immediate\write\gtoutfile{Title: \ifx\theasciititle\relax
\thetitle\else\theasciititle\fi}}
\immediate\write\gtoutfile{Subj-class: GT or SG, GR etc}
\immediate\write\gtoutfile{MSC-class: \theprimaryclass\ifx\thesecondaryclass\relax\else, \thesecondaryclass\fi}
\immediate\write\gtoutfile{Journal-ref: Algebr. Geom. Topol. \thevolumenumber\s
(\thevolumeyear) \startpage-\finishpage}
\immediate\write\gtoutfile{Comments: Published by Algebraic and
Geometric Topology at}
\immediate\write\gtoutfile{\s\s\s  http://www.maths.warwick.ac.uk/agt/AGTVol\thevolumenumber/agt-\thevolumenumber-\thepapernumber.abs.html}
\immediate\write\gtoutfile{\noexpand\\}
\immediate\write\gtoutfile{}
\ifx\theasciiabstract\relax
\immediate\write\gtoutfile{\theabstract}\else
\immediate\write\gtoutfile{\theasciiabstract}\fi
\immediate\write\gtoutfile{}
\immediate\write\gtoutfile{\noexpand\\}
\immediate\write\gtoutfile{}
\immediate\closeout\gtoutfile}}  

\def\maketitlepage{\makeagttitle\makeheadfile}
\let\makeshorttitle\maketitlepage
\let\maketitle\maketitlepage

\lognumber{44}
\volumenumber{2}
\volumeyear{2002}
\papernumber{44}
\pagenumbers{1147}{1154}
\received{19 July 2002}
\accepted{5 December 2002}
\published{19 December 2002}

\let\bigskip\medskip
\def\S{section~\ignorespaces}

\def\al{\alpha}

\def\hatg{\widehat{g}}
\def\ks{\widehat{KS}}
\def\al{\alpha}
\def\mn{M^n}
\def\btop{\rm BTOP}
\def\bstop{\rm BSTOP}

\def\bhml{\rm BHML}
\def\stop{\rm STOP}

\def\tgam{\widetilde{\gamma}}

\begin{document}

\title{Equivalences to the triangulation conjecture}

\author{Duane Randall}

\address{Department of Mathematics and Computer Science\\Loyola 
University, New Orleans, LA 70118, USA}

\email{randall@loyno.edu}

\begin{abstract}
We utilize the obstruction theory of
Galewski-Matumoto-Stern to derive equivalent formulations of the
Triangulation Conjecture. For example, every closed topological manifold
$\mn$ with $n\ge 5$ can be simplicially triangulated if and only if the
two distinct combinatorial triangulations of $RP^5$ are simplicially
concordant.
\end{abstract}
\asciiabstract{We utilize the obstruction theory of
Galewski-Matumoto-Stern to derive equivalent formulations of the
Triangulation Conjecture. For example, every closed topological manifold
M^n with n >4 can be simplicially triangulated if and only if the
two distinct combinatorial triangulations of RP^5 are simplicially
concordant.}


\maketitle

\section{Introduction}

The Triangulation Conjecture (TC) affirms that every closed topological
manifold $\mn$ of dimension $n\ge 5$ admits a simplicial triangulation.
The vanishing of the Kirby-Siebenmann class $KS(M)$ in $H^4(M;Z/2)$ is
both necessary and sufficient for the existence of a combinatorial
triangulation of $\mn$ for $n\ge 5$ by \cite{KS}. A combinatorial
triangulation of a closed manifold $\mn$ is a simplicial triangulation
for which the link of every $i$-simplex is a combinatorial sphere of
dimension $n-i-1$. Galewski and Stern \cite[Theorem 5]{GS1} and Matumoto
\cite{Ma} independently proved that a closed connected topological manifold
$\mn$ with $n\ge 5$ is simplicially triangulable if and only if
$$\delta_\al KS(M)=0~~~\hbox{in}~~~H^5(M;\ker\al)\leqno(1.1)$$
where $\delta_\al$ denotes the Bockstein operator associated to the exact
sequence\break $0\to\ker\al\to\theta_3\buildrel\al\over\longrightarrow
Z/2\to 0$ of abelian groups. Moreover, the Triangulation Conjecture is
true if and only if this exact sequence splits by \cite{GS1} or \cite[page 26]{Ra}.
The Rochlin invariant morphism $\al$ is defined on the homology bordism
group $\theta_3$ of oriented homology $3$-spheres modulo those which
bound acyclic compact $PL$ $4$-manifolds. Fintushel and Stern \cite{FS} and Furuta \cite{F}
proved that $\theta_3$ is infinitely generated.

We freely employ the notation and information given in Ranicki's excellent
exposition \cite{Ra}. The relative boundary version of the Galewski-Matumoto-Stern
obstruction theory in \cite{Ra} produces the following result. Given any
homeomorphism $f:|K|\to|L|$ of the polyhedra of closed $m$-dimensional $PL$
manifolds $K$ and $L$ with $m\ge 5$, $f$ is homotopic to a $PL$ homeomorphism if
and only if $KS(f)$ vanishes in $H^3(L;Z/2)$. More generally, a homeomorphism
$f:|K|\to|L|$ is homotopic to a $PL$ map $F:K\to L$ with acyclic point inverses
if and only if
$$\delta_\al(KS(f))=0~~~\hbox{in}~~H^4(L;\ker\al)~.\leqno(1.2)$$
Concordance classes of simplicial triangulations on $M^n$ for $n\ge 5$
correspond bijectively to vertical homotopy classes of liftings of the
stable topological tangent bundle $\tau:M\to\btop$ to $BH$ by
\cite[Theorem 1]{GS1} and so are enumerated by $H^4(M;\ker\al)$. The
classifying space $BH$ for the stable bundle theory associated to
combinatorial homology manifolds in \cite{Ra} is denoted by ${\rm
BTRI}$ in \cite{GS1} and by $\bhml$ in \cite{Ma}. We employ
obstruction theory to derive some known and new results and
generalizations of \cite{GS2} and \cite{S} on the existence of
simplicial triangulations in \S2 and to record some equivalent
formulations of $TC$ in \S3. Although some of these formulations may be
known, they do not seem to be documented in the literature.
\bigskip
\section{Simplicial Triangulations}
\bigskip
Let $\delta^*$ denote the integral Bockstein operator associated to
the exact sequence $0\to Z\buildrel \times 2\over\longrightarrow
Z\buildrel\rho\over\longrightarrow Z/2\to 0$. We proceed to derive
some consequences of the vanishing of $\delta^*$ on Kirby-Siebenmann
classes.  The coefficient group for cohomology is understood to be
$Z/2$ whenever omitted. Matumoto knew in \cite{Ma} that the vanishing
of $\delta^*KS(M)$ implied the vanishing of $\delta_\al KS(M)$. Let
$\iota_m$ denote the fundamental class of the Eilenberg-MacLane space
$K(Z,m)$. Since $H^{m+1}(K(Z,m);G)=0$ for all coefficient groups $G$,
trivially $\delta_\alpha(\rho\iota_m)=0$ in
$H^{m+1}(K(Z,m);\ker\alpha)$. Thus $\delta_\alpha$ vanishes on $KS(M)$
in (1.1) or $KS(f)$ in (1.2) whenever $\delta^*$ does. This
observation together with (1.1) and (1.2) justifies the following
well-known statements. Every closed connected topological manifold
$\mn$ with $n\ge 5$ and $\delta^*KS(M)=0$ admits a simplicial
triangulation. Let $f:|K|\to|L|$ be any homeomorphism of the polyhedra
of closed $m$-dimensional $PL$ manifolds $K$ and $L$ with $m\ge 5$. If
$\delta^*KS(f)=0$, then $f$ is homotopic to a $PL$ map $F:K\to L$ with
acyclic point inverses.

\bigskip
\noindent{\bf Proposition 2.1}\sl\qua All $k$-fold Cartesian products of closed
$4$-manifolds are simplicially triangulable for $k\ge 2$. All products $M^4\times S^1$ with non-orientable
closed $4$-manifolds $M^4$ are simplicially triangulable. Let $N^4$ be any simply connected closed $4$-manifold with $KS(N)$ trivial and also $b=\hbox{rank of}~H_2(N;Z)\ge 1$. Let $f:|K|\to |L|$ be any homeomorphism with $KS(f)$ nontrivial and $|K|=|L|=N\times S^1$. Then $f$ is homotopic to a $PL$ map $F:K\to L$ with acyclic point inverses.\rm

\bigskip
\noindent{\bf Proof of 2.1}\qua Since $KS(\gamma)$ is a primitive cohomology
class for the universal bundle $\gamma$ on $\btop$, we have $KS(M_1\times
M_2)=KS(M_1)\otimes 1+1\otimes KS(M_2)$ in $H^4(M_1\times M_2)$.
Triviality of $\delta_\al$ on $H^4(M^4)$ by dimensionality yields
triangulability of all $k$-fold products of closed $4$-manifolds for
$k\ge 2$, and of $M^4\times S^1$ by (1.1).

The product $N^4\times S^1$ admits $2^b$ distinct
combinatorial structures by \cite{KS}; moreover, for every non-zero class $u$
in $H^3(N\times S^1)$, there is a homeomorphism of polyhedra with
distinct combinatorial structures whose Casson-Sullivan invariant is $u$
by \cite[page 15]{Ra}. The vanishing of $\delta^*KS(f)$ follows from the
triviality of $\delta^*$ on $H^3(N\times S^1)=\rho(H^2(N;Z)\otimes
H^1(S^1;Z))$.\endproof
\medskip

No closed $4$-manifold $M^4$ with $KS(M)$ non-zero can be simplicially
triangulated. Yet $k$-fold products of such manifolds $M^4$ by (2.1) and
their products with spheres or tori produce infinitely many distinct
non-combinatorial, yet simplicially triangulable closed manifolds in
every dimension $\ge 5$. In contrast, there are no known examples of
non-smoothable closed $4$-manifolds which can be simplicially
triangulated, according to Problem 4.72 of \cite[page 287]{K}.

\bigskip
\noindent{\bf Theorem 2.2}\sl\qua Let $M^n$ be any closed connected
topological manifold with $n\ge 5$ such that the stable spherical
fibration determined by the tangent bundle $\tau(M)$ has odd order in
$[M,BSG]$. Suppose that either $H_2(M;Z)$ has no $2$-torsion or else
all $2$-torsion in $H_4(M;Z)$ has order $2$. Then $M$ is simplicially
triangulable.\rm

\bigskip
\noindent{\bf Proof}\qua The Stiefel-Whitney classes of $M$ are
trivial by the hypothesis of odd order. We first consider the special
case that $\tau(M)$ is stably fiber homotopically trivial. Let $g:M\to
SG/\stop$ be any lifting of a classifying map $\tau(M):M\to\bstop$ in
the fibration
$$
SG/\stop\buildrel j\over\longrightarrow
\bstop\buildrel\pi\over\longrightarrow BSG\leqno(2.3)
$$ 
The Postnikov $4$-stage of $SG/\stop$ is $K(Z/2,2)\times K(Z,4)$. Now
$j^*KS(\tgam)=\iota^2_2+\rho(\iota_4)$ by Theorem 15.1 of \cite[page
328]{KS} where $\tgam$ denotes the universal bundle over
$\bstop$. Clearly $\delta^*(j^*KS(\tgam))=\delta^*(\iota^2_2)=2u$
where $u$ generates $H^5(K(Z/2,2);Z)\approx Z/4$. If all
nonzero $2$-torsion in $H_4(M;Z)$ has order $2$, then
$\delta^*KS(M)=2g^*u=0$. If $H_2(M;Z)$ has no $2$-torsion, then
$\delta^*(g^*\iota_2)=0$ so again $\delta^*KS(M)=0$. Thus
$\delta_\al KS(M)=0$.

We suppose now that the stable spherical fibration of $\tau(M)$ has order $2a+1$ in $[M,BSG]$ with $a\ge 1$. Let $s:M\to S(2a\tau(M))$ be a section to the sphere bundle projection $p:S(2a\tau(M))\to M$ associated to $2a\tau(M)$. Now $S(2a\tau(M))$ is a stably fiber homotopically trivial manifold, since its stable tangent bundle is $(2a+1)p^*\tau(M)$. Since $KS(M)=(2a+1)KS(M)=s^*(KS(S(2a\tau(M))))$ we conclude that
$$\delta^*KS(M)=s^*(\delta^*KS(S(2a\tau(M))))=s^*0=0~.\leqno(2.4)$$
We consider the following homotopy commutative
diagram of principal fibrations.
$$\def\normalbaselines{\baselineskip20pt \lineskip3pt
\lineskiplimit3pt}
\def\mapright#1{\smash{\mathop{\longrightarrow}\limits^{#1}}}
\def\maprightlimits#1{\smash{\mathop{\longrightarrow}
\limits_{#1}}}
\def\mapdown#1{\Big\downarrow
\rlap{$\vcenter{\hbox{$\scriptstyle#1$}}$}}
\def\mapup#1{\Big\uparrow
\rlap{$\vcenter{\hbox{$\scriptstyle#1$}}$}}
\matrix{&&K(\ker\alpha,4)&\mapright{i}& (K(\ker\alpha,4),*)&=&\left(K(\ker\alpha,4),
*\right)\cr
&&{}\mapdown{}&&\mapdown{}&&\mapdown{i}\cr
&&BH&\mapright{i}&(BH,BPL)&\mapright{t}&\left(K(\theta_3,4),*\right)\cr
&&{}\mapdown{\pi}&&\mapdown{\hat\pi}&&\mapdown{\alpha}\cr
S^4&\mapright{ks}&BTOP&\mapright{i}&~(BTOP,BPL)
&\mapright{\widehat{KS}}&(K(Z/2,4),*)\cr
&&&&\mapdown{\delta_\alpha\widehat{KS}}
&&\mapdown{\delta_\alpha{\iota}}\cr
&&&&(K(\ker\alpha,5),*)
&=&(K(\ker\alpha,5),*)\cr}\leqno(2.5)$$
The fiber map $\alpha$ is induced
from the path-loop fibration on $K(\ker\alpha,5)$ via the Bockstein
operator $\delta_\alpha \iota$ on the fundamental class $\iota$ of
$K(Z/2,4)$. The induced morphism $\al_*$ on $\pi_4$ is the Rochlin
morphism $\al:\theta_3\to Z/2$ by construction. The relative principal
fibration $\hat\pi$ is induced from $\al$ via the map $\widehat{KS}$
classifying the relative universal Kirby-Siebenmann class. Thus
$(\widehat{KS}\circ i)^*\iota=KS(\gamma)$. Inclusion maps are denoted by
$i$ in (2.5). The induced morphisms $t_*$ and $(\widehat{KS})_*$ are
isomorphisms on $\pi_4$. We employ (2.5) in the proof of Theorem 3.1.

\section{Equivalent formulations to $TC$}

Galewski and Stern constructed a non-orientable closed connected
$5$-manifold $M^5$ in \cite{GS2} such that $Sq^1KS(M)$ generates
$H^5(M)\approx Z/2$. They also proved that any such $M^5$ is
``universal'' for $TC$. Moreover, Theorem 2.1 of \cite{GS2} essentially
affirms that either $TC$ is true or else no closed connected topological
$n$-manifold $M^n$ with $Sq^1KS(M)\not=0$ and $n\ge 5$ can be
simplicially triangulated.

\bigskip
\noindent{\bf Theorem 3.1}\sl

The following statements are equivalent to the
Triangulation Conjecture.
\begin{enumerate}
\item[\rm(1)] Any (equivalently all) of the classes $\delta_\al
KS(\gamma)$, $\delta_\al\ks$, and $\delta_\al \iota$ in (2.5) is
trivial if and only if any (equivalently all) of the fiber maps $\pi$,
$\hat\pi$, and $\al$ in (2.5) admits a section.

\item[\rm(2)]The essential map $f:S^4\cup_{2}e^5\to\btop$ lifts to $BH$ in (2.5).

\item[\rm(3)]$Sq^1KS(\hat\gamma)\not=0$ in $H^5(BH)$ for the universal bundle
${\hat\gamma}=\pi^*\gamma$ on $BH$.

\item[\rm(4)]Any closed connected topological manifold $M^n$ with $Sq^1KS(M)\not=0$
and $n\ge 5$ admits a simplicial triangulation.

\item[\rm(5)]Every homeomorphism $f:|K|\to |L|$ with $KS(f)$ non-trivial is
homotopic to a $PL$ map with acyclic point inverses where $K$ and $L$ are
any combinatorially distinct polyhedra with $|K|=|L|=N^4\times RP^2$.
Here $N^4$ denotes any simply connected, closed
$4$-manifold with $KS(N)$ trivial and positive rank for $H_2(N;Z)$.

\item[\rm(6)] All combinatorial triangulations of each closed connected $PL$
manifold $\mn$ with $n\ge 5$ are concordant as simplicial triangulations.

\item[\rm(7)] The two distinct combinatorial triangulations of $RP^5$ are
simplicially concordant.

\item[\rm(8)] Every closed connected topological manifold $M^n$ with $n\ge 5$ that is stably fiber homotopically trivial admits a simplicial triangulation.
\end{enumerate}\rm

\bigskip
\noindent{\bf Proof}\qua $TC\Leftrightarrow (1)$\qua Statement (1) is
equivalent to the splitting of the exact sequence $0\to\ker\al\to\theta_3\buildrel\al\over\longrightarrow
Z/2\to 0$ through the induced morphisms on homotopy in dimension
4.

\bigskip
$TC\Leftrightarrow (2)$\qua Let $ks:S^4\to\btop$
represent the Kirby-Siebenmann class in homotopy. That is, $[ks]$ has
order 2 and is dual to $KS(\gamma)$ under the $\bmod~2$ Hurewicz
morphism. Now $ks$ admits an extension $f:S^4\cup_2e^5\to\btop$, since
the cofibration exact sequence
$$\pi_5(\btop)\longrightarrow[S^4\cup_2 e^5, \btop]\to\pi_4(\btop)\buildrel
\times 2\over\longrightarrow
\pi_4(\btop)\leqno(3.2)$$
corresponds to $0\longrightarrow Z/2\longrightarrow Z\oplus Z/2\buildrel
\times 2\over\longrightarrow Z\oplus Z/2$. If $g:S^4\cup_2 e^5\to BH$ is any
lifting of $f$, the composite map using (2.5)
$$h:S^4\subset S^4\cup_{2}e^5\buildrel g\over\longrightarrow BH\buildrel
i\over\longrightarrow (BH, BPL)\buildrel t\over\longrightarrow
(K(\theta_3,4),*)\leqno(3.3)$$
produces $u=[h]$ in $\theta_3$ with $2u=0$ and $\al(u)=1$, since $\al(u)=[\al\circ
h]=[\ks\circ ks]$ generates $\pi_4(K(Z/2,4))$. Thus $TC$ is true. Conversely, if $TC$
is true, a section $s:\btop\to BH$ to $\pi$ in (2.5) gives a lifting $s\circ f$ of $f$.

\bigskip
$TC\Leftrightarrow (3)$\qua Properties of $KS(\gamma)$ are enumerated in
\cite{M} and \cite{R}. Since\break $Sq^1KS(\gamma)\not=0$, a section $s$ to
$\pi$ in (2.5) gives $Sq^1(KS(\hat\gamma)\not=0$ so $TC$ implies
$3$. We now assume that $TC$ is false and claim that the generator $Sq^1
\iota$ for $H^5(K(Z/2,4))\approx Z/2$ lies in the image of
$$H^5(K(\ker\al,5))\approx Hom(\pi_5(K(\ker_\al,5)),Z/2)\approx
Hom(\ker\al,Z/2).$$ The Serre exact sequence then gives
$\al^*(Sq^1\iota)=0$ in $H^5(K(\theta_3,4))$ so
$$Sq^1KS(\hat\gamma)=(t\circ i)^*(\al^* Sq^1 \iota)=0.$$ Thus we must
construct a morphism $\ker\al\to Z/2$ which does not extend to
$\theta_3$. We consider the sequence $\ker\al\buildrel
\times 2\over\longrightarrow
\ker\al\buildrel\rho\over\longrightarrow\ker\al\otimes Z/2$ and define
$h:\ker\al\otimes Z/2\to Z/2$ as follows. $h(v)=1$ if and only if
$v=\rho(2z)$ for some $z\in\theta_3$ with $\al(z)=1$. Now $h$ is a
well-defined and non-trivial morphism, since $\theta_3$ does not have an
element $u$ with $2u=0$ and $\al(u)=1$ by hypothesis. The composite
morphism $h\circ\rho:\ker\al\to Z/2$ does not extend to $\theta_3$.

\bigskip
$TC\Leftrightarrow (4)$\qua Suppose $M^n$ with $Sq^1 KS(M)\not=0$
admits a simplicial triangulation. Now $Sq^1
KS(M)=g^*Sq^1 KS(\hat{\gamma})$ for any lifting $g:M\to BH$ of
$\tau:M\to\btop$. Since $Sq^1 KS(\hat{\gamma})\not=0$, $TC$ holds
by (3).

\bigskip
$TC\Leftrightarrow (5)$\qua Clearly triviality of $\delta_\al\ks$ in
(2.5) gives $\delta_\al KS(f)=0$ via naturality for every $f$. Suppose
that $\delta_\al KS(f)=0$ for any such $f$ in 5. Now
$KS(f)=\rho(v)\otimes i^*a$ in $\rho(H^2(M;Z))\otimes H^1(RP^2)\approx
H^3(L)$. Here $a$ generates $H^*(RP^\infty)$ and $i:RP^2\subset
RP^\infty$. Naturality via the universal example $CP^\infty\times
RP^\infty$ for $\rho(v)\otimes i^*a$ gives $\delta_\al
KS(f)=v\otimes\delta_\al(i^* a)$. Since $i^*:H^2(RP^\infty;\ker\al)\to
H^2(RP^2;\ker\al)$ is a monomorphism, $\delta_\al(i^* a)=0$ if and
only if $\delta_\al(a)=0$. Now $\delta_\al(a)=0$ if and only if $TC$
is true via the fibration $$K(\ker\al, 1)\longrightarrow
K(\theta_3,1)\buildrel\al\over\longrightarrow RP^\infty.$$

\bigskip
$TC\Leftrightarrow (6)\Leftrightarrow (7)$\qua $TC$ holds if and only if
$\delta_\al\iota=0$ for the fundamental class $\iota$ of $K(Z/2,3)$.
Concordance classes of simplicial triangulations of $\mn$ arising from
combinatorial triangulations differ by classes in $\delta_\al H^3(M)$.
This subgroup of $H^4(M,\ker\al)$ is trivial by naturality if
$\delta_\al\iota=0$. Conversely, $\delta_\al H^3(RP^5)=0$ if the two
distinct combinatorial triangulations of $RP^5$ given by Theorem 16.5 in
\cite[pages 332 and 337]{KS} are simplicially concordant. But
$\delta_\al(a^3)=0$ if and only if $\delta_\al\iota=0$ via the skeletal
inclusion $RP_3^5\subset K(Z/2,3)$ and naturality for $RP^5\to RP^5_3$.
\bigskip

$TC\Leftrightarrow (8)$\qua Similar to Theorem 5.1 of \cite{Ru}, we
 consider a regular neighborhood of the 9-skeleton of $SG/\stop$
 embedded in $R^m$ for some $m\ge 19$ in order to obtain a smoothly
 parallelizable manifold $W$ with boundary and a map $g:W\to SG/\stop$
 which is a homotopy equivalence through dimension $7$. The double
 $DW$ is smoothly parallelizable and admits an extension
 $\widehat{g}:DW\to SG/\stop$.  Note that $(\hatg)^*$ is a
 monomorphism through dimension 7.  Let $h:M\to DW$ be a degree one
 normal map. Now $M$ is stably fiber homotopically trivial and $h^*$
 is a monomorphism in cohomology.  In particular, $(\hatg\circ h)^*$
 is a monomorphism on $H^5(SG/\stop;\ker\alpha)$.  We conclude that
 $\delta_\alpha KS(M)=(\hatg\circ h)^*(\delta_\alpha\iota^2_2)=0$ if
 and only if $\delta_\alpha\iota^2_2=0$ for the fundamental class
 $\iota_2$ of $K(Z/2,2)$. So statement (8) yields
 $\delta_\al\iota^2_2=0$.

Let $f:K(Z/2,2)\to K(Z/2,4)$ classify $\iota^2_2$. Since
$\delta_\al\iota^2_2=0$ assuming statement (8), $f$ admits a lifting
$h:K(Z/2,2)\to K(\theta_3,4)$ in (2.5) such that $f=\al\circ h$. The
diagram
$$\def\normalbaselines{\baselineskip20pt \lineskip3pt
\lineskiplimit3pt}
\def\mapright#1{\smash{\mathop{\longrightarrow}\limits^{#1}}}
\def\maprightlimits#1{\smash{\mathop{\longrightarrow}
\limits_{#1}}}
\def\mapdown#1{\Big\downarrow
\rlap{$\vcenter{\hbox{$\scriptstyle#1$}}$}}
\def\mapup#1{\Big\uparrow
\rlap{$\vcenter{\hbox{$\scriptstyle#1$}}$}}
\matrix{&&[CP^3, K(\theta_3,4)]&\approx&\theta_3&&\cr
&{~~~}\nearrow^{\kern-15pt{h_*}}{}~~&{}\mapdown{\alpha_*}&&\mapdown{}\cr
Z/2\approx[CP^3,K(Z/2,2)]&\mapright{f_*}&[CP^3,K(Z/2,4)]&\approx&Z/2
\cr}\leqno(3.4)$$
yields a splitting to the exact sequence $0\to\ker\alpha\to\theta_3\to Z/2\to 0$ so $TC$ holds.

\Addressesr
\end{document}

\bye